\documentclass[11pt,twoside]{amsart}
\usepackage{amssymb}

\newtheorem{thm}{Theorem}[section]
\newtheorem{prop}[thm]{Proposition}
\newtheorem{lem}[thm]{Lemma}
\newtheorem{cor}[thm]{Corollary}

\newtheorem{conj}[thm]{Conjecture}

\newcommand{\skipit}[1]{{}}
\newcommand{\prfend}{\hbox to7pt{\hfil}
\par\vskip-\baselineskip\hbox to\hsize
{\hfil\vbox {\hrule width6pt height6pt}}\vskip\baselineskip}

\newcommand{\myarrow}[2]{\hbox to #1pt{\hfil$\to$\hfil}{\hskip-#1pt{\raise
10pt\hbox to#1pt{\hfil$\scriptscriptstyle #2$\hfil}}}}

\begin{document}

\title{\Large Interval Conjectures \\for level Hilbert functions}

\author{\Large Fabrizio Zanello\\
{\ }\\
\small Department of Mathematics, University of Notre Dame\\
Notre Dame, IN 46556, USA\\E-mail: zanello@math.kth.se\\
Current address: Department of Mathematical Sciences\\
Michigan Technological University\\
Fisher Hall, Office 226C\\
Houghton, MI 49931-1295, USA.}

\maketitle

\markboth{F.\ Zanello}{Interval Conjectures for level Hilbert functions}

\begin{abstract}
We conjecture that the set of all Hilbert functions of (artinian) level algebras enjoys a very natural form of regularity, which we call the {\em Interval Conjecture} (IC): If, for some positive integer $\alpha $, $(1,h_1,...,h_i,...,h_e)$ and $(1,h_1,...,h_i+\alpha ,...,h_e)$ are both level $h$-vectors, then $(1,h_1,...,h_i+\beta ,...,h_e)$ is also level for each integer $\beta =0,1,..., \alpha .$ In the Gorenstein case, i.e. when $h_e=1$, we also supply the {\em Gorenstein Interval Conjecture} (GIC), which naturally generalizes the IC, and basically states that the same property simultaneously holds for any two symmetric entries, say $h_i$ and $h_{e-i}$, of a Gorenstein $h$-vector.

These conjectures are inspired by the research performed in this area over the last few years. A series of recent results seems to indicate that it will be nearly impossible to characterize explicitly the sets of all Gorenstein or of all level Hilbert functions. Therefore, our conjectures would at least provide the existence of a very strong - and natural - form of order in the structure of such important and complicated sets.

We are still far from proving the conjectures at this point. However, we will already solve a few interesting cases, especially when it comes to the IC, in this paper. Among them, that of Gorenstein $h$-vectors of socle degree 4, that of level $h$-vectors of socle degree 2, and that of non-unimodal level $h$-vectors of socle degree 3 and any given codimension.

\end{abstract}
{\ }\\

\section{Introduction and statement of the conjectures}

The goal of this paper is to present, and to start investigating, a conjecture - which we will call the {\em Interval Conjecture} (IC) - on the regularity of the structure of (artinian) level Hilbert functions (or $h$-vectors), as well as an extension thereof - the {\em Gorenstein Interval Conjecture} (GIC) - to the Gorenstein case. Indeed, Gorenstein and level algebras are at the center of a significant amount of research today in commutative algebra, also due to the applications they have to several other areas of mathematics, such as invariant theory, algebraic geometry, combinatorics, and even complexity theory.

Let us first introduce the notation we will be using in this note. We consider standard graded artinian algebras $A=R/I$, where $R=k[x_1,...,x_r]$, $I$ is a homogeneous ideal of $R$ and the $x_i$'s all have degree 1. We will assume that $k$ is a field of characteristic zero.

The {\em $h$-vector} of $A$ is $h(A)=h=(h_0,h_1,...,h_e)$, where $h_i=\dim_k A_i$ and $e$ is the last index such that $\dim_k A_e>0$. (Notice that $h$-vector and Hilbert function essentially coincide in the artinian case, so we may use either of the two denominations indifferently). Since we may suppose that $I$ does not contain non-zero forms of degree 1, $r=h_1$ is defined as the {\em codimension} of $A$.
 
The {\em socle} of $A$ is the annihilator of the maximal homogeneous ideal $\overline{m}=(\overline{x_1},...,\overline{x_r})\subseteq A$, namely $soc(A)=\lbrace a\in A {\ } \mid {\ } a\overline{m}=0\rbrace $. Since $soc(A)$ is a homogeneous ideal, we can define the {\em socle-vector} of $A$ as $s(A)=s=(s_0,s_1,...,s_e)$, where $s_i=\dim_k soc(A)_i$. Notice that $h_0=1$, $s_0=0$ and $s_e=h_e>0$. The integer $e$ is called the {\em socle degree} of $A$ (or of $h$). The {\em type} of the socle-vector $s$ (or of the algebra $A$) is type$(s)=\sum_{i=0}^es_i$.
 
If $s=(0,0,...,0,s_e=t)$, we say that the algebra $A$ (or its $h$-vector) is {\em level} (of type $t$). In particular, the case $t=1$ is called {\em Gorenstein}.

For the history of level algebras, which began with the seminal paper of Stanley \cite{St}, as well as for a very accurate account on the subject and its applications updated to the year 2003, we refer the reader to the memoir \cite{GHMS}. For applications to complexity theory (the connection being the theory of inverse systems, which will be introduced next), see for example \cite{NW} and \cite{KS}.

It was initially hoped that the Hilbert functions of level algebras had some predictable behavior, and that therefore they would have been eventually classified. Instead, even in the Gorenstein case, now the situation seems unclear and particularly complicated. The property that Stanley and Iarrobino initially asked whether it always held (it is indeed called the SI-condition), then turned out to be false in codimension $r=5$ and higher. Also, in type greater than one, even in codimension three, a complete list of level Hilbert functions seems really difficult to determine.

Recall that an $h$-vector $(h_0=1,h_1,h_2,...,h_e)$ is defined to be {\em an SI-sequence} if it is symmetric (i.e., $h_i=h_{e-i}$ for all $i$) and its first half is {\em differentiable} (i.e., the vector $(1,h_1-h_0,h_2-h_1,..., h_{\left \lfloor \frac{e}{2}\right \rfloor }-h_{\left \lfloor \frac{e}{2}\right \rfloor  -1})$ is also the $h$-vector of some artinian algebra, or, equivalently, satisfies Macaulay's theorem).

For Gorenstein Hilbert functions of codimension 2 and 3 a classification has been provided by Stanley: Indeed, they are exactly the SI-sequences starting with $(1,2,...)$ or $(1,3,...)$ (see \cite{St}; the simpler case $r=2$ was already known to Macaulay \cite{Ma}, and for $r=3$ Stanley's result was recently reproved by this author in an elementary fashion; see \cite{Za1}). Instead, for all codimensions $r\geq 5$, the SI-condition turned out not to hold; in particular, even unimodality (which is obviously a weaker condition than SI) was proven false for $r\geq 5$, first by Stanley and then by other authors (see \cite{St,BI,Bo,BL}). In codimension 4, the situation is more unclear, since  to date both SI and unimodality have neither been proved nor disproved, in spite of some progress (see \cite{IS} and \cite{MNZ2}).

The purpose of this paper is to present a conjecture of regularity for the set of level Hilbert functions (the IC), and its natural generalization (the GIC) to the Gorenstein case. Since it appears nearly impossible to characterize those sets explicitly, our conjectures would at least prove the existence of a very strong - and very natural - form of regularity within their structures.

\begin{conj}\label{conje}
i) (The Interval Conjecture (IC)). Suppose that, for some positive integer $\alpha $, both $(1,h_1,...,h_i,...,h_e)$ and $(1,h_1,...,h_i+\alpha ,...,h_e)$ are level $h$-vectors. Then $(1,h_1,...,h_i+\beta ,...,h_e)$ is also level for each integer $\beta =0,1,..., \alpha .$

ii) (The Gorenstein Interval Conjecture (GIC)). Suppose that, for some positive integer $\alpha $, both $(1,h_1,...,h_i,...,h_{e-i},...,h_{e-1},h_e=1)$ and $(1,h_1,...,h_i+\alpha ,...,h_{e-i}+\alpha,...,h_{e-1},h_e=1)$ are Gorenstein $h$-vectors. Then $(1,h_1,...,h_i+\beta ,..., h_{e-i}+\beta,...,h_{e-1},h_e=1)$ is also Gorenstein for each integer $\beta =0,1,..., \alpha .$

\end{conj}
\noindent{\bf Remark.} i) Notice that, because of the symmetry of Gorenstein $h$-vectors, the GIC clearly implies the IC for the Gorenstein case.

ii) The analogous conjecture to the IC in general does not hold for level algebras if more than one entry is changed at the time. For instance, both $h=(1,3,6,10,4)$ and $h^{'}=(1,3,6,9,3)$ are level $h$-vectors, but we cannot \lq \lq fill the gap" between $h$ and $h^{'}$, since $(1,3,6,10,3)$ is not level (these facts are easy to prove, for instance using inverse systems).

iii) Similarly, it is not possible to extend the GIC to more than one adjacent entry (modulo symmetry). For instance, $(1,3,5,7,7,5,3,1)$ and $(1,3,6,8,8,6,3,1)$ are both Gorenstein $h$-vectors, but $(1,3,5,8,8,5,3,1)$ is clearly not (it violates Macaulay's theorem).
\\
\\\indent
We will see in the next section that the SI-condition forces the form of regularity we conjecture in both the IC and the GIC. But even in those instances where the SI-condition fails, no \lq \lq irregularity" feature is apparently known. Notice that, because of the symmetry property, for Gorenstein $h$-vectors the only spot where the IC could possibly fail is the middle entry (and therefore only in case the socle degree is even).

This is one of the main reasons why, in the Gorenstein case, we also provide the more general GIC. Again, no counterexample to the GIC seems to have been provided or suggested in the literature so far.

For instance, Stanley \cite{St2} (see also the recent preprint \cite{MNZ1}) has proved the existence of the Gorenstein $h$-vector $(1,24,20,24,1)$, and it is easy to see that $(1,24,24,24,1)$ too is Gorenstein. A very natural \lq \lq IC-type" question therefore is: Is $(1,24,b,24,1)$ also a Gorenstein $h$-vector for each $b=21,22,23$? We will prove later in this paper that the answer is affirmative.

Determining all possible Hilbert functions of level algebras of type greater than 1 is even more difficult. Actually, a characterization is again known in codimension $r=2$ (see Iarrobino \cite{Ia1}), but in this case already in codimension 3 it seems hopeless to find an explicit characterization (see our recent paper \cite{Za4}, where we prove the existence of non-unimodal level Hilbert functions in any codimension $r\geq 3$). On the other hand, again, no irregularity property (in the same sense as above) has ever been proved, or even suggested by the several classes of level $h$-vectors discovered over the last few years (see, e.g., the Appendices of the memoir \cite{GHMS}, where the authors classify all level $h$-vectors of codimension 3 and small socle degree).

It seems too hard to prove the IC or the GIC in general at this point. However, there are several interesting cases that we are already able to settle, especially when it comes to the IC. In particular, we will prove it for Gorenstein $h$-vectors of socle degree 4, for level $h$-vectors of socle degree 2, and for non-unimodal level $h$-vectors of socle degree 3 and any given codimension. We will also begin to investigate the IC for level $h$-vectors with any socle degree.

The main tool we will use in this note is the theory of inverse systems (which will be defined in the next section), along with some technical results of commutative algebra, in particular of Iarrobino's (\cite{Ia1,Ia2}) and of ours (\cite{Za5}).

\section{The Gorenstein case}

In this section we will study the IC and the GIC in the case of Gorenstein $h$-vectors.

It is well-known that a Gorenstein $h$-vector $h$ is symmetric. In particular, it is not possible to change the entry of degree $i$ of $h$ and obtain another Gorenstein $h$-vector, without also changing the entry of degree $e-i$ of $h$ by the same value. This immediately implies a first, strong regularity result in the Gorenstein case: Namely, if an irregularity occurs in the sense of Conjecture \ref{conje}), i) (the IC), then it can only happen in degree $i=\frac{e}{2}$, i.e. exactly in the middle of $h$, and in particular its socle degree $e$ must be even. This is also one of the main reasons why we have introduced the GIC for Gorenstein $h$-vectors.

Let us now prove that the SI-condition, as we mentioned in the introduction, implies the IC and the GIC.

\begin{prop}\label{SI} Fix two positive integers $e$ and $r$, and suppose that the Gorenstein $h$-vectors of socle degree 
$e$ and codimension $r$ coincide with the SI-sequences of the form $(1,r,h_2,...,h_e)$. Then the GIC (and therefore the IC) holds for this class of Gorenstein $h$-vectors.
\end{prop}

{\bf Proof.} We first prove the IC. As we have just noticed, by symmetry $e$ must be even, and it suffices to show the result in degree $i=\frac{e}{2}$. Hence suppose that, for some positive integer $\alpha $, both $h^{(0)}:=(1,r,h_2,...,h_{\frac{e}{2}},...,h_e=1)$ and $h^{(\alpha )}:=(1,r,h_2,...,h_{\frac{e}{2}}+\alpha ,...,h_e=1)$ are Gorenstein $h$-vectors (i.e., SI-sequences, by hypothesis), and fix an integer $\beta $, $0\leq \beta \leq \alpha .$ We want to show  that $h^{(\beta )}=(1,r,h_2,...,h_{\frac{e}{2}}+\beta ,...,h_e=1)$ is also an SI-sequence.

By definition of SI, we have that
$$(\Delta h)^{(0)}:=\left(1,r-1,h_2-r,..., h_{\frac{e}{2}-1}-h_{\frac{e}{2}-2}, h_{\frac{e}{2}}-h_{\frac{e}{2}-1}\right )$$
and
$$(\Delta h)^{(\alpha )}:=\left(1,r-1,h_2-r,...,h_{\frac{e}{2}-1}-h_{\frac{e}{2}-2}, h_{\frac{e}{2}}+ \alpha -h_{\frac{e}{2}-1}\right)$$ are both $h$-vectors of artinian algebras. Hence, in their last degree, $\frac{e}{2}$, they must satisfy Macaulay's theorem, i.e. $h_{\frac{e}{2}}+ \alpha -h_{\frac{e}{2}-1}$ is bounded from above by some function of $h_{\frac{e}{2}-1}-h_{\frac{e}{2}-2}$ (we avoid recalling Macaulay's theorem precisely, since it is a well-known result and its technical statement is not needed here). Therefore, for $0\leq \beta \leq \alpha $, also $h_{\frac{e}{2}}+ \beta -h_{\frac{e}{2}-1}\leq h_{\frac{e}{2}}+ \alpha -h_{\frac{e}{2}-1}$ is bounded from above by the same function. This immediately implies that
$$(\Delta h)^{(\beta )}:=\left(1,r-1,h_2-r,...,h_{\frac{e}{2}-1}-h_{\frac{e}{2}-2}, h_{\frac{e}{2}}+ \beta -h_{\frac{e}{2}-1}\right)$$  satisfies Macaulay's theorem as well, and is therefore the $h$-vector of some artinian algebra. Hence $h^{(\beta )}$ is also an SI-sequence, as we wanted to show.

The above proof basically generalizes to prove the GIC {\em mutatis mutandis}, for any socle degree $e$. Indeed, by symmetry, it suffices to consider any degree $i\leq \left \lfloor \frac{e}{2}\right \rfloor $ in place of $\frac{e}{2}$, and the argument to prove the SI condition for $h^{(\beta )}$ in degrees less than or equal to $i$ does not change. The argument to prove the SI from degree $i$ to degree $i+1$ (in case $i< \left \lfloor \frac{e}{2}\right \rfloor $) is entirely similar and relies on the simple fact, concerning the first difference of $h^{(\beta )}$, that if a pair of entries, say $(a,b)$, satisfies Macaulay's theorem in degrees $i$ and $i+1$, then clearly the pair $(a+1,b-1)$ does too in the same degrees. This concludes the proof of the proposition.{\ }{\ }\qed
\\
\\\indent
In particular, by \cite{St}, we immediately have:

\begin{cor}\label{stst}
The IC and the GIC both hold for the class of Gorenstein $h$-vectors of codimension $r\leq 3$.
\end{cor}

The two main results of this section are contained in the following theorem:

\begin{thm}\label{gor} i) Suppose that
$$h=\left(1,h_1=r,h_2=\binom{r+1}{2},...,h_{\frac{e}{2}-1}=\binom{r+\frac{e}{2}-2}{\frac{e}{2}-1},h_{\frac{e}{2}}=a,h_{\frac{e}{2}+1}=h_{\frac{e}{2}-1},...,h_{e-1}=h_1,h_e=1 \right)$$
is a Gorenstein $h$-vector of codimension $r$ and even socle degree $e$, where all entries of $h$ up to degree $\frac{e}{2}-1$ are maximal (i.e., they coincide with the Hilbert function of a polynomial ring in $r$ variables). Then
$$\left(1,h_1,...,h_{\frac{e}{2}-1},b,h_{\frac{e}{2}+1},...,h_e\right)$$
is also Gorenstein for all integers $b$ such that $a\leq b\leq \binom {r+\frac{e}{2}-1}{\frac{e}{2}},$ the upper bound being sharp.

ii) Suppose that
$$h=\left(1,h_1=r,h_2=\binom{r+1}{2},...,h_{\left \lfloor \frac{e}{2}\right \rfloor -1}=\binom{r+\left \lfloor \frac{e}{2}\right \rfloor -2}{\left \lfloor \frac{e}{2}\right \rfloor -1},h_{\left \lfloor \frac{e}{2}\right \rfloor }=a,h_{\left \lfloor \frac{e}{2}\right \rfloor +1}=a,...,h_{e-1}=h_1,1 \right)$$
is a Gorenstein $h$-vector of codimension $r$ and odd socle degree $e$, where all entries of $h$ up to degree $\left \lfloor \frac{e}{2}\right \rfloor -1$ are maximal. Then
$$\left(1,h_1,...,h_{\left \lfloor \frac{e}{2}\right \rfloor -1},b,b,h_{\left \lfloor \frac{e}{2}\right \rfloor +2},...,h_{e-1},1\right)$$
is also Gorenstein for all integers $b$ such that $a\leq b\leq \binom {r+\left \lfloor \frac{e}{2}\right \rfloor -1}{\left \lfloor \frac{e}{2}\right \rfloor },$ the upper bound being sharp.
\end{thm}

Notice that the most interesting cases are when $h$ is not unimodal, since constructing the  $h$-vector of the statement for $b$ greater than or equal to the previous entry is a fairly standard exercise (for instance using inverse systems, which will be defined below).

In particular, Theorem \ref{gor}, i) completely settles the IC for Gorenstein $h$-vectors of socle degree 4:

\begin{cor}\label{corgor}
The IC holds for the class of
Gorenstein $h$-vectors of socle degree 4. Precisely, if $h=(1,r,a,r,1)$ is a Gorenstein $h$-vector, then $h=(1,r,b,r,1)$ is also Gorenstein for all integers $b$ such that $a\leq b\leq \binom{r+1}{2}$, the upper bound being sharp.
\end{cor}
\noindent
{\bf Remark.} In \cite{MNZ2}, Migliore, Nagel and this author solved a conjecture of Stanley (\cite{St2}) predicting the asymptotic behavior of the least value, $f(r)$, that the degree 2 entry of a socle degree 4 Gorenstein $h$-vector may assume. Namely, we proved that, for $r$ going to infinity, $f(r)\sim_r (6r)^{\frac{2}{3}}$.

Corollary \ref{corgor} now shows that {\em all} values can be attained by $h_2$ between the minimum possible, $f(r)$, and the maximum possible, $\binom{r+1}{2}$.

However, providing a characterization of all Gorenstein $h$-vectors of socle degree 4 - that in the light of our result has just become equivalent to determining $f(r)$ - still seems too hard a problem to solve for every given codimension $r$.
\\
\\
{\bf Example.} We know from \cite{St2} that $(1,40,30,40,1)$ is a Gorenstein $h$-vector.

Hence Corollary \ref{corgor} shows that $(1,40,b,40,1)$ is also Gorenstein for all $b=30,31,...,\binom{41}{2}=820$.
\\
\\\indent
Before proving Theorem \ref{gor}, we need to introduce {\em {(Macaulay's)} inverse systems}, also known as {\em Matlis duality}. For a complete introduction to the theory of inverse systems, we refer the reader to \cite{Ge} and \cite{IK}.
\\
\\
{\bf Inverse systems.} Given the polynomial ring $R=k[x_1,x_2,...,x_r]$, let $S=k[y_1,y_2,...,y_r]$, and consider $S$ as a graded $R$-module, where the action of $x_i$ on $S$ is partial differentiation with respect to $y_i$.
 
There is a one-to-one correspondence between artinian algebras $R/I$ and finitely generated $R$-submodules $M$ of $S$, where $I=Ann(M)$ is the annihilator of $M$ in $R$ and, conversely, $M=I^{-1}$ is the $R$-submodule of $S$ which is annihilated by $I$ (cf. \cite{Ge}, Remark 1), p.17).

If $R/I$ has socle-vector $s=(0,s_1,...,s_e)$, then $M$ is minimally generated by $s_i$ elements of degree $i$, for $i=1,2,...,e$, and the $h$-vector of $R/I$ is given by the number of linearly independent partial derivatives obtained in each degree by differentiating the generators of $M$ (cf. \cite{Ge}, Remark 2), p.17).
 
In particular, level algebras of type $t$ and socle degree $e$ correspond to $R$-submodules of $S$ minimally generated by $t$ elements of degree $e$.
\\
\\\indent
We now need two technical results, both due to Iarrobino (\cite{Ia1,Ia2}).

The first one (along with the lemma preceding it) determines the behavior of the $h$-vector of a given inverse system module $M$, once a form $F$, being the sum of the powers of $m$ general linear forms, is added to $M$. (Iarrobino's result actually holds under more general hypotheses than those we are going to state here.)

\begin{lem} (\cite{Ia1}, Proposition 4.7) \label{lem1}
Let $F=\sum_{p=1}^m L_p^e$ be a form of degree $e$ in $S=k[y_1,...,y_r]$, where the $L_p$'s are general linear forms (i.e., their coefficients are chosen in a non-empty Zariski-open subset of $k^{r}$). Then the Gorenstein algebra $R/Ann(F)$ has $h$-vector
$$h(m)=(1,h_1(m),...,h_e(m)=1),$$
where, for $j=1,...,e$,
$$h_j(m)=\min \lbrace m,\dim_kR_j,\dim_kR_{e-j}\rbrace .$$
\end{lem}
 
\begin{thm}(\cite{Ia1}, Theorem 4.8 A)\label{ia1}
Let $h=(1,h_1,...,h_e)$ be the $h$-vector of a level algebra $A=R/I$, where $I$ annihilates the $R$-submodule $M$ of $S$. Let $m\leq \binom{r-1+e}{e}-h_e$. Then, if $F$ is the sum of the $e$-th powers of $m$ general linear forms (the Gorenstein $h$-vector of $R/Ann(F)$, $h(m)$, is given by Lemma \ref{lem1}), then the level algebra associated to $M^{'}=\langle M,F \rangle $ has $h$-vector $H=(1,H_1,...,H_e),$ where, for $i=1,...,e$, $$H_i=\min \left \lbrace h_i+h_i(m),\binom{r-1+i}{i} \right \rbrace .$$
\end{thm}

The other result of Iarrobino that we will need is a lower bound on the entries of the $h$-vector of the general Gorenstein quotient of a given type 2 level algebra having the same socle degree.

\begin{thm}(\cite{Ia2}, Theorem 2.2)\label{ia2}
Let $A$ be a level algebra of type 2 and socle degree $e$ having $h$-vector $h=(1,h_1,h_2...,h_e)$, and let $H^{1,gen}=(1,H_1^{1,gen},H_2^{1,gen},...,H_e^{1,gen})$ be the $h$-vector of the general Gorenstein quotient of $A$ having socle degree $e$, say corresponding to an inverse system module $\langle F\rangle $. Let $G$ be another degree $e$ form of $M$ such that $M=\langle F,G\rangle $, and denote by $d_i$ the dimension of the space of common partial derivatives in degree $i$ of $F$ and $G$. Then, for each $i=1,2,...,e$,
$$H_i\geq h_{e-i}-d_i.$$
\end{thm}

We are now ready to prove Theorem \ref{gor}.
\\
\\\indent
{\bf Proof of Theorem \ref{gor}.} i) Let $D_1$ be a degree $e$ form of $S=k[y_1,...,y_r]$ which generates an inverse system module with the $h$-vector $h$ of the statement. For any given $b$, $a\leq b\leq \binom{r+\frac{e}{2}-1}{\frac{e}{2}}$, consider a form $D_2$ which is the sum of the $e$-th powers of $m=b-a$ general linear forms of $S$. Let $M=\langle D_1,D_2\rangle $. 

Then, in particular, by Theorem \ref{ia1} the $h$-vector, $h(M)$, of $M$ up to degree ${\frac{e}{2}-1}$ is equal to\\$\left(1,h_1=r,h_2=\binom{r+1}{2},...,h_{\frac{e}{2}-1}=\binom{r+\frac{e}{2}-2}{\frac{e}{2}-1}\right)$ (since $h$ was already maximal in those degrees), and in degree $\frac{e}{2}$ it is equal to $h(M)_{\frac{e}{2}}=a+(b-a)=b$. 

Also, in all degrees greater than or equal to $\frac{e}{2}$, the spaces spanned by the derivatives of $D_1$ and $D_2$ only have a trivial intersection. Hence, with the notation of Theorem \ref{ia2}, $d_i=0$ for all $i=\frac{e}{2},\frac{e}{2}+1,...,e$.

Therefore, by Theorem \ref{ia2}, for $i=\frac{e}{2},\frac{e}{2}+1,...,e$, a general Gorenstein quotient of $M$ of socle degree $e$ has its degree $i$ entry $H_i\geq h(M)_{e-i}$. Hence, by symmetry, $H_i= H_{e-i}=h(M)_{e-i}$, since clearly the $h$-vector of a quotient cannot be bigger than the $h$-vector of the original algebra. 

This suffices to prove that the general Gorenstein quotient of $M$ has the desired $h$-vector, $(1,h_1,...,h_{\frac{e}{2}-1},b,h_{\frac{e}{2}+1},...,h_e)$. The sharpness of the upper bound is obvious, since, in codimension $r$, $\binom {r+ \frac{e}{2} -1}{\frac{e}{2} }$ is clearly the highest possible entry that we can have in degree $\frac{e}{2}$.

ii) The argument for this case is entirely similar to the previous one, and therefore will be omitted.{\ }{\ }\qed

\section{The IC for level $h$-vectors with small socle degree}

This section is devoted to the study of the Interval Conjecture for level $h$-vectors of socle degree 2 and 3. Our first result, which solves the IC when $e=2$, uses a lemma of Bigatti and Geramita \cite{BG} on the minimum possible entries that a level $h$-vector may assume (see also our paper \cite{Za3}), and a theorem of Cho and Iarrobino \cite{CI}.

Let $n$ and $i$ be positive integers. Recall that the {\em i-binomial expansion of n} is
$$n_{(i)}:=\binom{n_i}{i}+\binom{n_{i-1}}{i-1}+...+\binom{n_j}{j},$$
where $n_i>n_{i-1}>...>n_j\geq j\geq 1$.

Under these hypotheses, the $i$-binomial expansion of $n$ is unique (e.g., see \cite{BH}, Lemma 4.2.6).

Following \cite{BG}, define, for any integer $a$,
$$\left(n_{(i)}\right)_{a}^{a}=\binom{n_i+a}{i+a}+\binom{n_{i-1}+a}{ i-1+a}+...+\binom{n_j+a}{ j+a}.$$

The following two results are phrased in a different way from the original papers.

\begin{lem} (\cite{BG}, Lemma 3.3) \label{bg}
Let $h_d$ be the degree $d$ entry of an $h$-vector $h$. Then the least possible value that $h$ may assume in degree $d-1$ is
$$\left(\left(h_d\right)_{(d)}\right)_{-1}^{-1}.$$
\end{lem}

\begin{prop} (\cite{CI}, Theorem 1.4) \label{ci}
Let $M\subset S=k[y_1,...,y_r]$ be the inverse system module corresponding to a level algebra of socle degree $d$. Then the degree $d-1$  entry of the $h$-vector of $M$, $h_{d-1}$, satisfies the following inequalities:
$$\left(\left(h_d\right)_{(d)}\right)_{-1}^{-1}\leq h_{d-1}\leq \min \left \lbrace \binom{r+d-2}{d-1},r\cdot h_d \right \rbrace .$$
Moreover, each integral value of $h_{d-1}$ within the above range is attained by such a module $M$.  
\end{prop}

We are now ready to prove the following:

\begin{prop}\label{2}
{\ }\\
The IC holds for the class of level $h$-vectors of socle degree 2. Precisely:\\
i) For every given positive integer $t$, the $h$-vector $(1,r,t)$ is level for all integers $r\geq \left(t_{(2)}\right)_{-1}^{-1}$, the lower bound being sharp.
\\\noindent
ii) For every given positive integer $r$, the $h$-vector $(1,r,t)$ is level for all integers $t$ such that $1\leq t\leq \binom{r+1}{2}$, the upper bound being sharp.
\end{prop}

{\bf Proof.} i) The sharpness of the lower bound for $r$ follows from Lemma \ref{bg}. That all values of $r$ as in the statement are attained is a simple consequence of Proposition \ref{ci} and the fact that we work in any codimension.

ii) For any codimension $r$, it suffices to consider a level inverse system generated by $t$ forms, each being the sum of the squares of $r$ general linear forms, and the result follows from Theorem \ref{ia1}.{\ }{\ }\qed
\\
\\
{\bf Example.} Let $t=25$. Then $(1,r,25)$ is a level $h$-vector for all integers $r\geq 7$.

Indeed, the 2-binomial expansion of 25 is $25=\binom {7}{2}+\binom{4}{1}$. Hence
$$\left(25_{(2)}\right)^{-1}_{-1}=\binom {6}{1}+\binom{3}{0}=6+1=7.$$
\\\indent
The next theorem concerns the much more complicated case of level $h$-vectors of socle degree 3. The \lq \lq IC-type" result that we will prove on the second entry of such $h$-vectors includes, in particular, the non-unimodal case.

Let us first introduce a recent theorem of ours, which supplies a lower bound for the $h$-vector of the general type $c$ level quotient of a given level algebra of type $t$ having the same socle degree.

\begin{thm}(\cite{Za5}, Theorem 2.9)\label{za}
Let $A$ be a type $t$ level algebra with $h$-vector $h=(1,h_1,h_2,...,h_e)$, and let, for $c=1,2,...,t-1$, $H^{c,gen}:=(1,H_1^{c,gen},H_2^{c,gen},...,H_e^{c,gen})$ be the $h$-vector of the general type $c$ level quotient of $A$ having the same socle degree $e$. Then, for each $i=1,2,...,e $, we have:
$$H_i^{c,gen}\geq \frac{1}{ t^2-1}\left((t-c)h_{e-i}+(ct-1)h_i\right).$$
\end{thm}

We are now ready to show:

\begin{thm}\label{3}
Let $(1,r,a,t)$ be a level $h$-vector, such that
$$t(r-2t)+3\leq a\leq \binom{r+1}{2}-1 {\ }\mbox{ and }{\ } r\geq t(a-t)+2.$$
Then $(1,r,a+1,t)$ is also level.
\end{thm}

{\bf Proof.} {\em Claim.} It is enough to show that, under the hypotheses of the statement, $h=(1,r_1,a,t-1)$ is level whenever $r_1\in \lbrace r-1,r\rbrace $.

{\em Proof of Claim.} If $r_1=r$, by Theorem \ref{ia2} it suffices to add, to an inverse system module $M$ with $h$-vector $h$, the third power of one general linear form, and by Theorem \ref{ia1} we obtain the level $h$-vector $(1,r,a+1,t)$, as desired.

If $r_1=r-1$, it suffices to add to $M$ the third power of a new variable, and by inverse systems it is easy to see that we again obtain $(1,r,a+1,t)$. This proves the Claim.

Therefore, using Theorem \ref{za} applied to the general type $c=t-1$ level quotient of a level algebra with $h$-vector $(1,r,a,t)$, since the entries of an $h$-vector must be integers, it is enough to prove that the inequalities:
$$\frac{1}{t^2-1}(a+r(t(t-1)-1))>r-2$$
and
$$\frac{1}{t^2-1}(r+a(t(t-1)-1))>a-1$$
are both satisfied.

But a standard computation shows that these two inequalities are equivalent to the lower bounds for $a$ and $r$ of the hypotheses of the statement, and the theorem follows.{\ }{\ }\qed
\\
\\\indent
One of the most interesting applications of the above result is to non-unimodal $h$-vectors. Indeed, by iterating Theorem \ref{3}, we immediately have:

\begin{cor}\label{cor3}
Let $(1,r,a,t)$ be a level $h$-vector, where $t\geq \frac{r}{2}$. Then $(1,r,b,t)$ is also level for all integers $b=a,a+1,...,\min \lbrace t+1,\binom{r+1}{2}\rbrace .$
\end{cor}

The range for $b$ in the above result can be extended, provided that we suppose that $t$ be closer to $r$. Indeed:

\begin{thm}\label{33}
Let $(1,r,a,t)$ be a level $h$-vector, where $t\geq r-2$. Then $(1,r,b,t)$ is also level for all integers $b=a,a+1,...,\min \lbrace r\cdot t,\binom{r+1}{2}\rbrace ,$ the upper bound being sharp.
\end{thm}

{\bf Proof.} We can suppose that $r-2\geq \frac{r}{2}$, since otherwise the few cases not already covered by Corollary \ref{cor3} are trivial. Hence the case $b\leq t+1$ has already been proved in Corollary \ref{cor3}.

Instead, for any $b$, $b\geq t+2\geq r$, let us consider $t$ forms of degree 3, being the sum of powers of, respectively, $m_1$,..., $m_t$ general linear forms, where $m_1+m_2+...+m_t=b$, and $1\leq m_i\leq r$ for all $i=1,2,...,t$. By inverse systems and Theorem \ref{ia2}, the result now easily follows.

The upper bound of the statement for $b$ is clearly sharp, since $t$ forms in $r$ variables can have a total of at most $r\cdot t$  first partial derivatives. This completes the proof.{\ }{\ }\qed
\\
\\\noindent
{\bf Example.} It can be shown, using Stanley's construction provided in \cite{St2}, that $(1,40,30,40,1)$ is  a Gorenstein $h$-vector. Thus, $(1,40,30,40)$ is level (in general, it is a standard fact that the truncation of a level $h$-vector is again a level $h$-vector - a very simple proof can be given, e.g.,  by inverse systems).

By using sums of powers of linear forms, we also obtain, for instance, that $(1,40,35,45)$ is level. Hence, by Theorem \ref{33}, we have that $(1,40,b,45)$ is a level $h$-vector for all integers  $b=35,36,...,\\
\min \lbrace \binom{41}{2}, 40\cdot 45 \rbrace =820.$

\section{The IC for level $h$-vectors with any socle degree}

The purpose of this section is to study the IC for the last entry of level $h$-vectors of any socle degree. 
We first show that part ii) of Proposition \ref{2} easily extends to any socle degree $e$. (As usual, $\lceil x\rceil $ will denote the smallest integer grater than or equal to $x$.)

\begin{prop}\label{e}
Given positive integers $r$ and $e$, the $h$-vector
$$h=\left(1,h_1=r,h_2=\binom{r+1}{2},...,h_{e-1}= \binom{r+e-2}{e-1},h_e=t\right )$$
is level for each integer $t$ such that $\left \lceil \frac{\binom{r+e-2}{e-1}}{r} \right \rceil \leq t\leq \binom{r+e-1}{e}$, both bounds being sharp.
\end{prop}

{\bf Proof.} The upper bound is trivial, and the lower bound is a simple consequence, by inverse systems, of the fact that any form in $r$ variables can have at most $r$ first partial derivatives.

By Theorem \ref{ia1}, if for each integer $t$ within the range of the statement, we consider a level module generated by $t$ forms, each being the sum of the $e$-th powers of $r$ general linear forms, the proposition easily follows.{\ }{\ }\qed
\\
\\\indent
The following is the main result of this section. A fundamental tool will be again Theorem \ref{za}. We have:

\begin{thm}\label{ee}
Let $h=(1,h_1,...,h_{e-1},h_e=t)$ be a level $h$-vector, and suppose that, for each $i=1,2,...,e-1$, the following inequality is satisfied:
$$h_{e-i}+t^2-th_i-1>0.$$
Then $(1,h_1,...,h_{e-1},t-1)$ is also a level $h$-vector.
\end{thm}

{\bf Proof.} Let us apply Theorem \ref{za} to the $h$-vector of the type $c=t-1$ general level quotient of socle degree $e$ of a level algebra having $h$-vector $h$. Hence it clearly suffices to prove, for each  $i=1,2,...,e-1$, that
$$\frac{1}{t^2-1}(h_{e-i}+h_i(t(t-1)-1))>h_i-1.$$
But it is easy to check that these inequalities are equivalent to those of the statement, and the proof is complete.{\ }{\ }\qed
\\
\\\indent
An interesting particular case of the previous result is the following:

\begin{cor}\label{eecor}
Let $h=(1,h_1,...,h_{e-1},h_e=t_0)$ be a level $h$-vector, where $t_0\geq M:=\max \lbrace h_1,h_2,\\...,h_{e-1}\rbrace $. Then $(1,h_1,...,h_{e-1},t)$ is also level for each $t=M-1,M,...,t_0$.
\end{cor}

{\bf Proof.} We may assume that $h_j\geq 2$ for each $j=1,2,...,e-1$, otherwise clearly $h$ cannot be level for $r>1$.
Hence  the result follows by iterating the application of Theorem \ref{ee}, since its hypotheses are clearly satisfied for each value of $t$ between $M$ and $t_0$.{\ }{\ }\qed
\\
\\\noindent
{\bf Example.} Using the construction of \cite{St2}, one obtains the Gorenstein $h$-vector $(1,60,45,40,45,60,1)$. Hence, by truncation, $(1,60,45,40,45,60)$ is level. If we add, for instance, 10 fifth powers of one general linear form each, we have that $h=(1,60,55,50,55,70)$ is also level.

Therefore, by Corollary \ref{eecor} applied to $h$, it follows that $(1,60,55,50,55,t)$ is a level $h$-vector for all integers $t=59,60,...,70$.\\
\\
{\bf Acknowledgements.} We wish to thank Juan Migliore and Uwe Nagel for dozens of insightful and inspiring conversations on Hilbert functions and level algebras that we have had during the last year. We also thank the referee for very helpful comments.


\end{document}